\documentclass{amsart}
\usepackage{amsfonts,amssymb,amsmath,amsthm}
\usepackage{url}
\usepackage[dvips]{epsfig}
\newtheorem{thrm}{Theorem}[section]
\newtheorem{lem}[thrm]{Lemma}
\newtheorem{prop}[thrm]{Proposition}

\theoremstyle{definition}
\newtheorem{definition}[thrm]{Definition}

\author{Carlo~Alberto~Mantica and Luca~Guido~Molinari}
\address{Physics Department\\
Universit\'a degli Studi di Milano and I.N.F.N. sez. Milano\\
Via Celoria 16, 20133 Milano, Italy.\\
Present address of C.~A.~Mantica: I.I.S. Lagrange, Via L. Modignani 65, 
20161, Milano, Italy}
\email{carloalberto.mantica@libero.it, luca.molinari@mi.infn.it}
\subjclass[2010]{Primary 53B20, Secondary 53B21}
\keywords{Codazzi tensor, Riemann tensor, generalized curvature tensor, 
recurrent 1-forms.}
\begin{document}
\title{Extended Derdzinski-Shen theorem\\ for the Riemann tensor}
\begin{abstract} We extend a classical result by Derdzinski and Shen, on 
the restrictions imposed on the Riemann tensor by the existence of a 
nontrivial Codazzi tensor. The new conditions of the 
theorem include Codazzi tensors (i.e. closed 1-forms) 
as well as tensors with gauged Codazzi condition (i.e. ``recurrent 1-forms''), 
typical of some well known differential structures. 
\end{abstract}
\maketitle
\section{Introduction}
Codazzi tensors are of great interest in the geometric literature and 
have been studied by several authors, as Berger and Ebin \cite{[1]}, 
Bourguignon \cite{[3]}, Derdzinski \cite{[5],[6]}, Derdzinski and Shen 
\cite{[7]}, Ferus \cite{[8]}, Simon \cite{[14]}; a compendium of results 
is reported in Besse's book \cite{[2]}.\\ 
In the following, $M$ is a $n\ge 4$ dimensional Riemannian manifold with 
metric $g_{ij}$ and Riemannian connection $\nabla $; the Ricci tensor is
$R_{kl}=-R_{mkl}{}^m$ and the scalar curvature is $R=g^{ij}R_{ij}$ \cite{[19]}. 
%
A $(0,2)$ symmetric tensor is a {\em Codazzi tensor} if it satisfies 
the Codazzi equation:
\begin{equation}
\nabla_j b_{kl} - \nabla_k b_{jl}=0.\label{Codazzi}
\end{equation}
A Codazzi tensor is {\em trivial} if it is a constant multiple of the 
metric tensor \cite{[7]}. 
%
In terms of differential forms, the Codazzi equation is the condition
for closedness of the 1-form $B_j = b_{jk} dx^k$, with covariant 
exterior differential $DB_j = \nabla_i b_{jk} \; dx^i\wedge dx^k$
\cite{[11],[3]}.

Codazzi tensors occur naturally in the study of harmonic Riemannian manifolds.
For example, the Ricci tensor is a Codazzi tensor    
if and only if $\nabla_m R_{jkl}{}^m=0$ (by the contracted second Bianchi 
identity), i.e. the manifold has harmonic Riemann curvature \cite{[2]}. 
The other way, the Weyl 1-form 
$\Sigma_j =\left ( R_{kj} -\frac{R}{2(n-1)} g_{kj}\right )dx^k$
is closed if and only if 
$\nabla_m C_{jkl}{}^m=0$, where $C_{jkl}{}^m$ is the conformal 
curvature tensor \cite{[13]}, 
i.e. the manifold has harmonic conformal curvature \cite{[2]}.

Berger and Ebin \cite{[1]} proved that in a compact Riemannian 
manifold the Codazzi tensors are parallel if the sectional curvature is 
non-negative and if there is a point where it is positive. In ref.\cite{[3]} 
the important 
geometric and topological consequences of the existence of a non-trivial 
Codazzi tensor are examined, particularly the restrictions 
imposed on the structure of the curvature operator. Derdzinki and Shen 
improved these results and proved a remarkable theorem on the consequences 
of the existence of a 
non trivial Codazzi tensor \cite{[7]}. The theorem is reported also in 
Besse's book \cite{[2]}. 

\begin{thrm}[Derdzinski-Shen]\label{thrm1}
Let $b_{ij}$ be a Codazzi tensor on a Riemannian manifold $M$, $x$ a point
of $M$, $\lambda $ and $\mu $ two eigenvalues of the operator $b_i{}^j(x)$,
with eigenspaces $V_{\lambda}$ and $V_\mu$ in $T_xM$. Then, the subspace 
$V_\lambda \wedge V_\mu$ 
is invariant under the action of the curvature operator $R_x$. 
\end{thrm}
The theorem can be rephrased as follows: given eigenvalues $\lambda$, 
$\mu$, $\nu$ of  $b_i{}^j(x)$ and corresponding eigenvectors $X$, $Y$, $Z$ 
in $T_xM$, it is $ R(X,Y)Z=0 $, provided that $\lambda$ and $\mu$ are different 
from $\nu$.\\
 
We point out that the Codazzi equation is a {\em sufficient} condition for
the theorem to hold. A more general condition is suggested by the lemma:
\begin{lem}\label{lem2.1}
Any $(0,2)$ symmetric Codazzi tensor $b_{kl}$ satisfies the 
algebraic identity:
\begin{equation}
b_{im} R_{jkl}{}^m+b_{jm} R_{kil}{}^m+b_{km} R_{ijl}{}^m=0.\label{eq2.1}
\end{equation}
\begin{proof}
The following condition involving commutators (indices are cyclically
permuted) is true for a Codazzi tensor:
\begin{equation}
[\nabla_j,\nabla_k]b_{il}+[\nabla_k,\nabla_i]b_{lj}+[\nabla_i,\nabla_j]b_{kl}=0
\label{commutators}
\end{equation}
Each commutator is evaluated: 
$[\nabla_i,\nabla_j]b_{kl}=R_{ijk}{}^mb_{ml}+R_{ijl}{}^mb_{km}$. 
Three terms cancel by the first Bianchi identity, and the result is obtained.
\end{proof}
\end{lem}
We note in passing that also the following equation is solved by Codazzi
tensors:
\begin{equation}
[\nabla_i,\nabla_j]b_{kl}+[\nabla_j,\nabla_k]b_{li}+[\nabla_k,\nabla_l]b_{ij}+
[\nabla_l,\nabla_i]b_{jk} =0
\end{equation} 
By evaluating the commutators and simplifying the result by means of the
first Bianchi identity and the Codazzi property, the equation gives:
\begin{equation}
R_{kij}{}^mb_{ml}+R_{jli}{}^mb_{mk}+R_{ljk}{}^mb_{mi}+R_{ikl}{}^mb_{mj}=0
\end{equation}

We are ready to formulate the main theorem of the paper: it states that if 
a symmetric tensor $b_{kl}$ satisfies the algebraic identity (\ref{eq2.1}), 
then the same conclusions of the Derdzinski-Shen theorem are valid. Moreover,
the proof is much simpler.

\section{An extension of the Derdzinski-Shen theorem.}
\begin{definition}\label{def2.1}
A $(0,4)$ tensor $K_{ijlm}$ is a {\em generalized curvature tensor} \cite{[10]}
if it has the symmetries of the Riemann curvature tensor:\\
a) $K_{ijkl} =-K_{jikl} =-K_{ijlk}$,\\ 
b) $K_{ijkl} =K_{klij}$,\\
c) $K_{ijkl} +K_{jkil} +K_{kijl} =0$ (first Bianchi identity).
\end{definition}
\begin{lem}\label{lem2.2}
If a symmetric tensor $b_{kl}$ satisfies eq.(\ref{eq2.1}), 
then $K_{ijkl} = R_{ijrs} b_k{}^r b_l{}^s$ is a generalized curvature tensor.
\begin{proof}
Properties a) are shown easily. For example:\\ 
$K_{ijlk}= R_{ijrs}b_l{}^r b_k{}^s = R_{ijsr} b_l{}^s b_k{}^r=
- R_{ijrs}b_l{}^s b_k{}^r =-K_{ijkl}$.\\ 
Property c) follows from the condition (\ref{eq2.1}): 
$ K_{ijkl} +K_{jkil} +K_{kijl} = 
R_{ijrs}b_k{}^r b_l{}^s + R_{jkrs}  b_i{}^r b_l{}^s + R_{kirs} b_j{}^r b_l{}^s
=  ( R_{jis}{}^r b_{kr}+ R_{kjs}{}^r b_{ir}+ R_{iks}{}^r b_{jr}) b_l{}^s=0 $.\\
Property b) follows from c): $K_{ijkl} +K_{jkil} +K_{kijl} =0$. Sum the
identity over cyclic permutations of all indices $i,j,k,l$ and use the 
symmetries a) (this fact was pointed out in \cite{[10]}).\\ 
It is easy to see that a first Bianchi identity holds also for the last 
three indices: $K_{ijkl} +K_{iklj} +K_{iljk} =0$.
\end{proof}
\end{lem}

\begin{thrm}[main theorem]\label{thrm2.1}
Let $M$ be a $n$-dimensional Riemannian manifold. If a $(0,2)$ symmetric 
tensor $b_{kl}$ satisfies the algebraic equation 
$$ b_{im} R_{jkl}{}^m+b_{jm} R_{kil}{}^m+b_{km} R_{ijl}{}^m=0 $$
then, for any $x\in M$, arbitrary eigenvalues $\lambda$ and $\mu$ of 
$b_i{}^j(x)$, with eigenspaces $V_\lambda$ and $V_\mu $ in $T_xM$, 
the subspace $V_\lambda\wedge V_\mu$ is invariant under the action of
the curvature operator $R_x$.
\end{thrm}
The theorem can be rephrased as follows:\\ 
{\em If $b_{rs}$ is a symmetric tensor with the property} (\ref{eq2.1})
{\em and $X$, $Y$ and $Z$ are three eigenvectors of the matrix $b_r{}^s$
at a point $x$ of the manifold, with eigenvalues $\lambda $, $\mu $ and $\nu $,
then}
\begin{equation}
X^i \, Y^j\, Z^k \, R_{ijkl} =0 \label{eq2.10}
\end{equation}
{\em provided that $\lambda $ and $\mu $ are different from $\nu$.} 
\begin{proof}
Consider the first Bianchi identity for the Riemann tensor, the condition 
eq.(\ref{eq2.1}) and the third Bianchi identity for the curvature 
$K_{\ell ijk}=R_{\ell irs}b_j{}^r b_k{}^s$, 
and apply them to the three eigenvectors. The resulting algebraic relations can 
be put in matrix form:
\begin{equation}
\left[ \begin{array}{ccc}
1 & 1 & 1 \\
\lambda & \mu & \nu \\
\mu\nu & \lambda\nu & \lambda\mu 
\end{array}\right]
\left [ \begin{array}{c}
R_{\ell ijk}X^iY^jZ^k\\
R_{\ell jki}X^iY^jZ^k\\
R_{\ell kij}X^iY^jZ^z
\end{array}\right]=
\left [ \begin{array}{c} 0\\ 0\\ 0 \end{array}\right]
\end{equation}
The determinant of the matrix is
$(\lambda-\mu)(\lambda-\nu)(\nu-\mu)$. If the eigenvalues are all different 
then $R_{\ell ijk}X^iY^jZ^k=0$; by the symmetries of the Riemann tensor the 
statement is true for the contraction of any three indices.\\
Suppose now that $\lambda=\mu\neq \nu$, i.e. $X$ and $Y$ belong to the same
eigenspace; the system of equations implies that $R_{\ell kij}X^iY^jZ^k=0$. 
\end{proof}

\begin{prop}
The hypothesis (\ref{eq2.1}) of the main theorem is fulfilled if the symmetric 
tensor $b_{ij}$ is a Codazzi tensor or, more generally, if it solves a 
{\em gauged Codazzi equation} 
\begin{equation}
(\nabla_k -\beta_k)b_{ij} = (\nabla_i-\beta_i)b_{kj}\label{gauge}
\end{equation}
where the covariant gauge field $\beta_k$ is closed:
$\nabla_k\beta_j-\nabla_j\beta_k=0$. 
\begin{proof}
The proof is straightforward and may start from eq.(\ref{commutators})
by noting that the closedness condition for the gauge field ensures that
$[\nabla_i -\beta_i,\nabla_j-\beta_j] = [\nabla_i, \nabla_j]$. Then
eq.(\ref{commutators}) is again true and eq.(\ref{eq2.1}) follows. 
\end{proof}
\end{prop}

The gauged Codazzi equation (\ref{gauge}) can be interpreted 
in terms of differential forms.

\section{Recurrent tensors forms}
\begin{definition}\label{def3.1}
A 1-form $B_j =b_{kj} dx^k$ is {\em recurrent} if there is a nonzero 
scalar 1-form $\beta =\beta_i dx^i$ such that
\begin{equation}
D B_j =\beta \wedge B_j.\label{eq3.2}
\end{equation}
\end{definition}
In local components it is
 $(\nabla_i -\beta_i)b_{kl} \,(dx^i\wedge dx^k)=0 $. Therefore, 
the 1-form $B_j =b_{kj} dx^k$ is recurrent if and only if \cite{[11]}
\begin{equation}
(\nabla_i-\beta_i) b_{kl} = (\nabla_k -\beta_k)b_{il}. \label{eq3.4}
\end{equation}
%
If $\beta =0$ the closedness of the 1-form $B_j =b_{kj} dx^k$ and  
the Codazzi equation are recovered. 
Eq.(\ref{eq3.4}) enlarges the ordinary notion of recurrence, by which a tensor
is recurrent if $\nabla_i b_{kl} =\beta_i b_{kl}$.
An example of ordinary recurrence are the Ricci recurrent spaces, where 
$\nabla_i R_{kl} =\beta_i R_{kl} $ \cite{[9]}.\\ 
The extended recurrence eq.(\ref{eq3.4}) includes well known differential 
structures, as the {\em Weakly b symmetric} manifolds, defined by the condition
\begin{equation}
\nabla_i b_{kl} = A_i b_{kl} + B_k b_{il} +D_l b_{ik} .\label{eq3.12}
\end{equation}
For these manifolds, eq.(\ref{eq3.4}) is 
satisfied with the choice $\beta_i =A_i -B_i$ and, if the covector 
$\beta_i $ is a closed 1-form, the main theorem applies. Weakly Ricci 
symmetric manifolds (see \cite{[4]} and \cite{[12]} 
for a compendium) are of this sort:
\begin{equation}
\nabla_i R_{kl} =A_i R_{kl} +B_k R_{il} +D_l R_{ik}.
\end{equation}
They were introduced by Tamassy and Binh \cite{[18]}.


\begin{thebibliography}{99}
%
\bibitem{[1]} M.~Berger and D.~Ebin, 
{\em Some characterizations of the space of symmetric tensors on a Riemannian manifold}, 
J. Differential Geometry \textbf{3} (1969) 379-392.
%
\bibitem{[2]} A.~L.~Besse, {\em Einstein Manifolds}, Springer (1987) .
%
\bibitem{[3]} J.~P.~Bourguignon, 
{\em Les vari\'et\'es de dimension 4 \`a signature non nulle dont la 
courbure est harmonique sont d'Einstein}, 
Invent. Math. \textbf{63} n.2 (1981) 263-286.
%
\bibitem{[4]} U.~C.~De, 
{\em On weakly symmetric structures on Riemannian manifolds}, 
Facta Universitatis, Series; Mechanics, Automatic Control and 
Robotics vol \textbf{3} n.14  (2003) 805-819.

\bibitem{[5]} A.~Derdzinski,
{\em Some remarks on the local structure of Codazzi tensors}, 
in Global differential geometry and global analysis, lecture notes 
n.\textbf{838}, Springer-Verlag (1981) 251-255.

\bibitem{[6]} A.~Derdzinski, 
{\em On compact Riemannian manifolds with harmonic curvature}, 
Math. Ann. {\bf 259} (1982) 145-152.

\bibitem{[7]} A.~Derdzinski and C.~L.~Shen, 
{\em Codazzi tensor fields, curvature and Pontryagin forms}, 
Proc. London Math. Soc. \textbf{47} n.3 (1983) 15-26.

\bibitem{[8]} D.~Ferus,
{\em A remark on Codazzi tensors in constant curvature spaces}, 
in global differential geometry and global analysis, lecture notes 
n.\textbf{838}, Springer-Verlag (1981) 257.

\bibitem{[9]} Q.~Khan, 
{\em On Recurrent Riemannian manifolds}, 
Kyungpook Math. J. \textbf{44} (2004) 269-276.

\bibitem{[10]} S.~Kobayashi and K.~Nomizu, 
{\em Foundations of differential geometry}, Vol. 1, Interscience (1963) 
New York.

\bibitem{[11]} D.~Lovelock and H.~Rund, 
{\em Tensors, differential forms and variational principles}, 
reprint Dover Ed. (1988).

\bibitem{[12]} C.~A.~Mantica and L.~G.~Molinari, 
{\em A second order identity for the Riemann tensor and applications}, 
Colloq. Math. {\bf 122} n.1 (2011) 69-82; arXiv:0802.0650 [math.DG].

\bibitem{[13]} M.~M.~Postnikov, 
{\em Geometry VI , Riemannian geometry}, Encyclopaedia of Mathematical 
Sciences, Vol.{\bf 91} (2001) Springer.

\bibitem{[14]} U.~Simon,
{\em Codazzi tensors}, in global differential geometry and global analysis, 
lecture notes n.\textbf{838}, Springer-Verlag (1981), 289-296 .

\bibitem{[18]} L.~Tamassy and T.~Q.~Binh, 
{\em On weakly symmetries of Einstein and Sasakian manifolds},
Tensor (N.S) \textbf{53} (1993) 140-148.

\bibitem{[19]} R.~C.~Wrede, 
{\em Introduction to vector and tensor analysis}, Reprint Dover (1963).
\end{thebibliography}
\end{document}